\RequirePackage{fix-cm}
\documentclass[smallcondensed,numbook]{svjour3}     
\smartqed  
\usepackage{graphicx}
\usepackage{url}
\usepackage{xcolor}
\definecolor{newcolor}{rgb}{.8,.349,.1}

\usepackage[hidelinks, colorlinks, allcolors=blue, bookmarksdepth=2]{hyperref}

\definecolor{gold}{rgb}{0.85, 0.65, 0.13}
\global\long\def\red#1{\textcolor{red}{#1}}%
\usepackage[numbered]{bookmark}

\usepackage{framed,multirow}

\usepackage{amsmath}
\usepackage{amssymb}
\usepackage{amsfonts}

\usepackage{mathrsfs}
\usepackage{upgreek}

\usepackage{epstopdf}
\usepackage{algorithmic}
\ifpdf
\DeclareGraphicsExtensions{.eps,.pdf,.png,.jpg}
\else
\DeclareGraphicsExtensions{.eps}
\fi
\usepackage[export]{adjustbox}

\usepackage{xspace}
\usepackage{bm}
\usepackage{physics}
\usepackage{cleveref}
\usepackage{enumerate}
\usepackage[caption=false]{subfig}
\usepackage{mathtools}
\usepackage{anyfontsize}

\usepackage{booktabs}
\usepackage{multirow, multicol, makecell}
\usepackage{longtable}
\crefname{assumption}{Assumption}{Assumptions}
\crefname{definition}{Definition}{Definitions}
\crefname{lemma}{Lemma}{Lemmas}
\crefname{remark}{Remark}{Remarks}
\crefname{theorem}{Theorem}{Theorems}
\crefname{proposition}{Proposition}{Propositions}
\crefname{section}{Section}{Sections}
\crefname{figure}{Fig.}{Figs.}
\crefname{equation}{}{}
\crefname{table}{Table}{Tables}
\crefname{appendix}{Appendix}{Appendices}

\newcommand{\fnc}[1]{\ensuremath{\mathcal{#1}}}

\renewcommand{\H}[0]{\mathsf{H}}

\newcommand{\E}[0]{\mathsf{E}}

\newcommand{\T}[0]{\mathsf{T}}

\newcommand{\A}[0]{\mathsf{A}}

\newcommand{\J}[0]{\mathsf{J}}

\newcommand{\Dxi}[0]{{\mathsf{D}}_{\hat{\bm{x}}_{i}}}

\newcommand{\Qxi}[0]{{\mathsf{Q}}_{\hat{\bm{x}}_{i}}}

\newcommand{\Exi}[0]{{\mathsf{E}}_{\hat{\bm{x}}_{i}}}

\newcommand{\Nxig}[0]{{\mathsf{N}}_{\hat{\bm{x}}_{i}\gamma}}

\newcommand{\R}[0]{\mathsf{R}}

\newcommand{\Bg}[0]{\mathsf{B}_{\gamma}}

\newcommand{\Rg}[0]{\mathsf{R}_{\gamma}}

\newcommand{\V}[0]{\mathsf{V}}

\newcommand{\etal}[0]{{et~al.\@}\xspace}
\newcommand{\eg}[0]{{e.g.\@}\xspace}
\newcommand{\ie}[0]{{i.e.\@}\xspace}

\newcommand{\ignore}[1]{} 

\newcommand{\polyref}[1]{\ensuremath{\mathbb{P}^{#1}}(\hat{\Omega})}

\newcommand{\IRtwo}[2]{\mathbb{R}^{{#1}\times{#2}}}%

\newcommand{\pder[2]}[]{\dfrac{\partial#1}{\partial#2}}%
\global\long\def\fn#1{\mathcal{#1}}%
\global\long\def\mathds#1{\mathds{#1}}%

\usepackage{etoolbox}
\newcommand*{\affaddr}[1]{#1} 
\newcommand*{\affmark}[1][*]{\textsuperscript{#1}}
\usepackage{tablefootnote}
\usepackage{threeparttable}
\usepackage{tablefootnote}

\usepackage{graphicx,scalerel}
\newcommand\sbullet[1][.7]{\mathbin{\ThisStyle{\vcenter{\hbox{%
					\scalebox{#1}{$\SavedStyle\bullet$}}}}}%
}
\usepackage[left=25mm,right=25mm,top=25mm,bottom=25mm]{geometry}

\begin{document}

\title{\vskip -0.5cm Quadrature Rules on Triangles and Tetrahedra for Multidimensional Summation-By-Parts Operators} 


\titlerunning{Quadrature Rules on Triangles and Tetrahedra}        

\author{Zelalem Arega Worku\protect\affmark[$\dagger$] \and Jason E. Hicken\protect\affmark[$\ddagger$] \and David W. Zingg\protect\affmark[$\dagger$]}

\authorrunning{Z. Worku, J.E. Hicken, D.W. Zingg}

\institute{Zelalem Arega Worku \at
	\email{zelalem.worku@mail.utoronto.ca}  
	\and
	Jason E. Hicken \at
	\email{hickej2@rpi.edu}
	\and
	David W. Zingg \at
	\email{dwz@oddjob.utias.utoronto.ca}
	\and
	\affaddr{\affmark[$\dagger$]Institute for Aerospace Studies, University of Toronto, Toronto, Ontario, M3H 5T6, Canada\\}	
	\affaddr{\affmark[$\ddagger$]Department of Mechanical, Aerospace, and Nuclear Engineering, Rensselaer Polytechnic Institute, Troy, NY, United States}	
}

\date{}
\def\makeheadbox{\relax}
\maketitle
\vskip -2.5cm
\addcontentsline{toc}{section}{Abstract}
\begin{abstract}
Multidimensional diagonal-norm summation-by-parts (SBP) operators with collocated volume and facet nodes, known as diagonal-$ \E $ operators, are attractive for entropy-stable discretizations from an efficiency standpoint. However, there is a limited number of such operators, and those currently in existence often have a relatively high node count for a given polynomial order due to a scarcity of suitable quadrature rules. We present several new symmetric positive-weight quadrature rules on triangles and tetrahedra that are suitable for construction of diagonal-$ \E $ SBP operators. For triangles, quadrature rules of degree one through twenty with facet nodes that correspond to the Legendre-Gauss-Lobatto (LGL) and Legendre-Gauss (LG) quadrature rules are derived. For tetrahedra, quadrature rules of degree one through ten are presented along with the corresponding facet quadrature rules. All of the quadrature rules are provided in a supplementary data repository. The quadrature rules are used to construct novel SBP diagonal-$ \E $ operators, whose accuracy and maximum timestep restrictions are studied numerically.


\end{abstract}

\section{Introduction} \label{sec:introduction}
Summation-by-parts (SBP) operators enable the construction of entropy-stable high-order discretizations of the Euler and Navier-Stokes equations \cite{fisher2013high,crean2018entropy,chan2018discretely,gassner2018br1,parsani2015entropy,fernandez2019staggered,ranocha2018comparison}.  Diagonal-norm SBP operators have collocated solution and volume quadrature nodes, which enables trivial inversion of the norm/mass matrix, facilitating efficient implementation of high-order methods with explicit time-marching schemes. Due to the collocation, the efficiency of the method is significantly affected by the number of quadrature nodes. This is even more prominent for the Hadamard-form entropy-stable discretizations \cite{fisher2013high,crean2018entropy} of the Euler and Navier-Stokes equations on simplices, as a volume flux computation coupling each node with all other nodes must be calculated. Therefore, the development of quadrature rules with fewer nodes on simplices is imperative to improve efficiency of entropy-stable discretizations with diagonal-norm SBP operators. While existing positive interior (PI) quadrature rules offer the fewest nodes for a given quadrature accuracy, their use for entropy-stable SBP discretizations requires expensive element coupling computations, although there are ways to reduce this cost to some extent \cite{chan2018discretely,shadpey2020entropy}. Alternatives to PI rules for SBP operators have been proposed; notably, Hicken \etal \cite{hicken2016multidimensional} derived quadrature rules with a set of volume nodes on each facet of the triangle and tetrahedron, and Chen and Shu \cite{chen2017entropy} presented rules on the triangle that enable construction of SBP operators with collocated volume and facet quadrature nodes, which are referred to as diagonal-$ \E $ or $ \R^{0}$\cite{marchildon2020optimization} SBP operators. Diagonal-$ \E $ SBP operators eliminate the need to extrapolate the solution from volume to facet nodes and are of particular interest for entropy-stable SBP discretizations as they reduce the cost of element coupling operations and enable straightforward enforcement of boundary conditions. However, the number of nodes required for the quadrature rules of existing diagonal-$ \E $ operators is significantly larger than that of the PI rules, especially for the tetrahedron. Furthermore, only a limited number of rules of this type are available in the literature. In light of this, the goal of this paper is to find efficient quadrature rules on the triangle and tetrahedron for the purpose of constructing efficient diagonal-norm diagonal-$ \E $ multidimensional SBP operators.

Quadrature rules with boundary nodes on simplices have been explored, although to a lesser extent than PI rules. Sets of nodes on the triangle that are typically well-suited for interpolation or quadrature accuracy have both been utilized to derive such rules, \eg, see \cite{taylor2000algorithm,cohen2001higher,mulder2001higher,giraldo2006diagonal,wingate2008performance,witherden2014analysis,liu2017higher} among others. Similar studies for the tetrahedron, however, are lacking. A shared property of several of the nodal sets in the mentioned studies is that for a degree $ p $ operator there are $ p+1 $ nodes on each facet of the triangle, and vertices are included. This automatically excludes the rules from being applicable for construction of diagonal-$ \E $ SBP operators, as the facet quadrature accuracy is not sufficient. We recall that a sufficient condition on the facet quadrature rule for the existence of a degree $ p $ SBP operator is that it be at least of degree $ 2p $ accurate. Although this is only a sufficient condition, as is evident from the existence of Legendre-Gauss-Lobatto (LGL) tensor-product operators on quadrilaterals and hexahedra, to the authors' knowledge, no diagonal-$ \E $ SBP operator on simplices with facet quadrature rule of degree less than $ 2p $ has been constructed.

Odd-degree quadrature rules on the triangle for diagonal-$ \E $ SBP operators with Legendre-Gauss (LG) facet nodes were first introduced in \cite{chen2017entropy}, followed by even degree rules with LGL and LG facet nodes in \cite{fernandez2019staggered,hicken2023summationbyparts} and odd-degree rules with LGL facet nodes in \cite{wu2021high}. Together, these studies provide quadrature rules up to degree eight\footnote{\cite{hicken2023summationbyparts} provides additional rules of degree 12 and 16 with LGL facet nodes.} on the triangle. For the tetrahedron, Hicken \cite{hicken2023summationbyparts} derived even degree quadrature rules up to degree eight with PI rule facet nodes. Marchildon and Zingg \cite{marchildon2020optimization} studied optimization of these operators and developed novel rules up to degree four on the tetrahedron, managing to lower the number of nodes for the degree two and four quadrature rules. This was achieved by allowing the facet nodes to be placed at the vertices and edges of the tetrahedron. Other efforts to improve the efficiency of entropy-stable discretizations with multidimensional SBP operators include, for instance, the use of staggered grids in  \cite{fernandez2019staggered}, entropy-split formulations in \cite{worku2023entropy}, and collapsed coordinate tensor-product elements in  \cite{montoya2023efficient}.

In this paper, we derive symmetric quadrature rules for construction of diagonal-norm diagonal-$ \E $ SBP operators on triangles and tetrahedra using the open-source Julia code \texttt{SummationByParts.jl }\cite{hicken2023summationbyparts}, which employs the Levenberg-Marquardt algorithm (LMA) \cite{levenberg1944method,marquardt1963algorithm} to solve the nonlinear systems of equations that arise from the quadrature accuracy conditions. The code's capability is enhanced by enforcing a constraint to find positive weights and by combining it with a particle swarm optimization (PSO) \cite{kennedy1995particle} subroutine to mitigate issues related to initial guesses and convergence to suboptimal local minima. We extend the available set of quadrature rules for diagonal-$ \E $ operators up to degree twenty on triangles with both the LGL and LG facet node configurations, and up to degree ten for tetrahedra, achieving a significant reduction in the number of nodes relative to many of the existing rules. The new rules are used to construct novel diagonal-$ \E $ SBP operators, whose accuracy and timestep restrictions are studied numerically.

The rest of the paper is organized as follows: \cref{sec:preliminaries} describes the problem statement and symmetry groups on the reference triangle and tetrahedron, \cref{sec:methodology} details the methodology employed, \cref{sec:quad rules} presents the derived quadrature rules along with a description of multidimensional SBP operators and their construction, and numerical results are presented in \cref{sec:num results} followed by conclusions in \cref{sec:conclusions}.

\section{Preliminaries}\label{sec:preliminaries}
Constructing quadrature rules over a domain of interest requires solving highly nonlinear systems of equations to find the nodal locations and weights. Usually, quadrature rules are designed to be exact for a desired degree of polynomial functions. The problem can be stated as: find $ \bm{x} $ and $ \bm{w} $ such that 
\begin{equation}\label{eq:prob}
	\int_{\Omega}\fn P_{j}\left(\widehat{\bm{x}}\right)\dd{\Omega}=\sum_{i=1}^{n_p}\bm{w}_{i}\fn P_{j}\left(\bm{x}_{i}\right),\qquad j\in\{1,\dots,n_{b}\},
\end{equation}
where $ \widehat{\bm{x}} = [x_{1},\dots,x_{d}]^{T} $, $ d $ is the spatial dimension, and  $ n_{p} $ and $ n_{b} $ denote the number of quadrature points and polynomial basis functions, respectively. For a degree $ q_{v} $ accurate quadrature rule on a simplex, there are $ n_{b}= {{p+d} \choose{d}} $ polynomial basis functions. Rewriting the problem statement in matrix form, we have
\begin{equation}\label{eq:prob matrix}
	 \bm{g} \coloneqq \V^{T}\bm{w}-\bm{f}=\bm{0},
\end{equation}
where $ \V $ is the Vandermonde matrix containing evaluations of the basis functions at each node along its columns and $ \bm{f}=\left[\int_{\Omega}\fn P_{1}\left(\widehat{\bm{x}}\right)\dd{\Omega},\dots,\int_{\Omega}\fn P_{n_{b}}\left(\widehat{\bm{x}}\right)\dd{\Omega}\right]^{T} $. The basis functions used to construct the Vandermonde matrix affect its condition number. It is well-known that, for high-order polynomials, monomial basis functions result in an ill-conditioned $ \V $. In contrast, the orthonormal Proriol-Koornwinder-Dubiner (PKD) \cite{proriol1957family,koornwinder1975two,dubiner1991spectral} basis functions offer better conditioning and a convenient $ \bm{f} $ vector; all except the first entry of $ \bm{f} $ are zero due to the orthogonality of the basis functions. 

For the purpose of constructing a degree $ p $ diagonal-norm diagonal-$ \E $ multidimensional SBP operator, we require that:
\begin{enumerate}[i.]
	\item all quadrature points lie in the closure of the simplex,
	\item all weights are positive,
	\item the volume quadrature is at least degree $ q_v=2p-1 $ accurate, 
	\item a subset of the quadrature points lying on each facet form a positive-weight facet quadrature rule of at least degree $q_{f}=2p $, and
	\item both the facet and volume quadrature rules are symmetric. 
\end{enumerate}
The symmetry requirement ensures that a solution obtained using the SBP operators is not spatially biased within an element.  Furthermore, the symmetry constraints reduce the number of unknowns in \cref{eq:prob matrix} substantially \cite{witherden2015identification}.

The reference triangle and tetrahedron are defined, respectively, as 
\begin{align} \label{eq:ref elems}
	 \Omega_{\text{tri}}&=\{\left(x,y\right)\mid x,y\ge-1; \; x+y\le0\},\\
	  \Omega_{\text{tet}}&=\{\left(x,y,z\right)\mid x,y,z\ge-1; \; x+y+z\le-1\}.
\end{align}
There are three symmetry groups on the triangle, and five on the tetrahedron \cite{felippa2004compendium,zhang2009set,witherden2015identification}. However, we will identify symmetric nodes that lie on the facets of the simplices as being in separate symmetry groups. On the triangle, the symmetry groups, in barycentric coordinates, are permutations of
\begin{equation}\label{eq:sym tri}
	\begin{aligned}
		S_{1} &= \left(\frac{1}{3}, \frac{1}{3}, \frac{1}{3}\right), &&
		S_{21} = (\alpha,\alpha,1-2\alpha),&&
		S_{111}  = (\alpha,\beta,1-\alpha-\beta),\\
		S_{\text{vert}} &= \left(1,0,0\right),&&
		S_{\text{mid-edge}}  = \left(\frac{1}{2},\frac{1}{2},0\right), && 
		S_{\text{edge}} =\left(\alpha,1-\alpha,0\right),
\end{aligned}
\end{equation}
where $ \alpha $ and $ \beta $ are parameters such that the quadrature points lie in the closure of the domain. The symmetry groups in the first line of \cref{eq:sym tri} represent interior points, while those in the second line denote points on the facets. Similarly, the symmetry groups on the reference tetrahedron are permutations of
\begin{equation}\label{eq:sym tet}
	\begin{aligned}
		&S_{1} = \left(\frac{1}{4},\frac{1}{4},\frac{1}{4},\frac{1}{4}\right), & &S_{\text{face-cent}} = \left(\frac{1}{3}, \frac{1}{3}, \frac{1}{3}, 0\right),\\
		&S_{31} = \left(\alpha,\alpha,\alpha,1-3\alpha\right), & &S_{\text{vert}}  = \left(1,0,0,0\right),\\
		&S_{22} = \left(\alpha,\alpha,1-\alpha,1-\alpha\right), & &S_{\text{mid-edge}}= \left(\frac{1}{2},\frac{1}{2},0, 0\right),\\
		&S_{211} = \left(\alpha,\alpha,\beta,1-2\alpha-\beta\right), & &S_{\text{face-}21}= (\alpha,\alpha,1-2\alpha, 0),\\
		&S_{1111} = \left(\alpha,\beta,\gamma,1-\alpha-\beta-\gamma\right), & &S_{\text{edge}} = \left(\beta,1-\beta, 0, 0\right),\\
		& & &S_{\text{face-}111} = \left(\alpha,\beta,1-\alpha-\beta, 0\right).
	\end{aligned}
\end{equation}
The facet symmetry groups in the right columns of \cref{eq:sym tet} are special cases of the symmetry groups in the left column, which are used in this work to denote interior points exclusively. 

\section{Methodology}\label{sec:methodology}
The Vandermonde matrix in \cref{eq:prob matrix} is a function of the quadrature points, $ \bm{x} $; hence, the algorithm to solve the equation starts by guessing the nodal locations and weights. This is done indirectly by providing the type and number of symmetry groups and the values of the associated parameters and weights. Using the initial guess of the parameters, it is possible to compute the coordinates of the $ i$-th node using the transformation,
\begin{equation} \label{eq:x=T lambda}
	\bm{x}_{i} = \T^{T} \bm{\lambda}_{k},
\end{equation}
where $ \T \in \IRtwo{(d+1)}{d}$ contains the coordinates of the $ d+1 $ vertices in its rows and $ \bm{\lambda}_{k} $ is the $ k $-th permutation of the barycentric coordinates of the symmetry group that corresponds to the $ i $-th node. The weight vector, $ \bm{w} $, is constructed by assigning equal weights to all nodes in the same symmetry group. 

To derive a degree $ q_{v} $ quadrature rule satisfying all the properties required to construct a degree $ p $ SBP diagonal-$ \E $ operator, we first need to find a facet quadrature rule of degree $ q_{f} \ge 2p $. On the reference triangle, we use either the LGL rule with $ n_{f} = p+2 $ facet nodes (including the vertices) or the LG rule with $ n_{f} = p+1 $. In construction of the volume quadrature rule, we fix the facet quadrature points; hence, the parameters in the facet symmetry groups are kept constant, \ie, we solve for the weights at all points and for the parameters associated with the interior symmetry groups. A similar strategy is followed to find the quadrature rules on the tetrahedron. While existing PI rules can be used as facet quadrature rules for the tetrahedron, they generally lead to more volume nodes than necessary. Hence, we first construct facet quadrature rules of degree $ 2p $ that would result in fewer volume quadrature points on the tetrahedron by placing some of the nodes at the vertices and/or edges of the facets. We note that, depending on their nodal locations, symmetry groups with the same number of nodes on the triangle can produce a different number of nodes on the tetrahedron when applied to its facets. For example, each of the $ S_{\text{vert}}  $ and $ S_{\text{face-}21} $ symmetry groups results in three nodes per facet, but four and twelve volume nodes, respectively. The number of volume nodes due to inclusion of the various facet symmetry groups of the tetrahedron is presented in \cref{tab:node contribution tet facet}. 

\begin{table*} [!t]
	\small
	\caption{\label{tab:node contribution tet facet} Number of tetrahedron volume node due to inclusion of facet symmetry groups}
	\centering
	\setlength{\tabcolsep}{1.1em}
	\renewcommand*{\arraystretch}{1.2}
	\begin{tabular}{ccccccc}
		\toprule
		group &  $ S_{\text{face-cent}} $&  $ S_{\text{vert}}  $ &  $ S_{\text{mid-edge}}$ & $ S_{\text{face-}21} $ & $ S_{\text{edge}}  $ & $ S_{\text{face-}111} $\\ 
		\midrule
		\# nodes & 4 & 4 & 6 & 12 & 12 & 24\\
		\bottomrule
	\end{tabular}
\end{table*}

The root finding problem in \cref{eq:prob matrix} can be written equivalently as a least-square minimization problem,
\begin{equation}\label{eq:min prob}
	\min_{\bm{\tau}}{\frac{1}{2}\bm{g}^T\bm{g}},
\end{equation}
where $ \bm{\tau} = [\widehat{\bm{\lambda}},\widehat{\bm{w}}]^T $, and $\widehat{\bm{\lambda}} $ and $\widehat{\bm{w}}$  are vectors of all the parameters and weights associated with each symmetry group.  The LMA is widely used to solve \cref{eq:min prob} as it is less sensitive to initial guesses than Newton's method. The LMA computes the step direction, $ \bm{h} $, which is initialized as the zero vector, as  
\begin{equation}\label{eq:step}
	\widetilde{\bm{h}} =-\widetilde{\A}^{+} \widetilde{\J^T \bm{g}},	
\end{equation}
where $ \A = \J^T \J + \nu \text{diag}(\J^T\J )$, $ \nu > 0$ controls the scale of exploration. Note that the notation $ (\cdot)^{+} $ in \cref{eq:step} denotes the Moore-Penrose pseudo-inverse, $\widetilde{(\cdot)} $ denotes the extraction of the rows and columns of a vector or matrix that correspond to parameters of the interior symmetry groups and all weights, $ \J\in \IRtwo{n_{b}}{n_{\tau}} $ is the Jacobian matrix given by
\begin{align}
	\J_{(i,j)} &= \pder[\bm{g}_{i}]{\bm{\tau}_{j}},
\end{align}
and $ n_{\tau} $ is the sum of the number of parameters and weights. The Jacobian can also be written in terms of matrices as
\begin{align}
	\J=\left[\sum_{k=1}^{d}\V_{x_{k}}^{T}\text{diag}(\bm{w})\pder[\bm{x}_{(:,k)}]{\widehat{\bm{\lambda}}},\V^{T}\pder[\bm{w}]{\widehat{\bm{w}}}\right],
\end{align}
where $ \V_{{x}_k} $ is the $ k $-th direction derivative of $ \V $ and $  \bm{x}_{(:,k)}  $is the $ k $-th direction component vector of $ \bm{x} $. The matrix $ \partial{\bm{x}_{(:,k)}}/\partial{\widehat{\bm{\lambda}}} $ is computed using the relation in \cref{eq:x=T lambda}, and $ \partial{\bm{w}}/\partial{\widehat{\bm{w}}} $ is a matrix of zeros and ones. The value of $ \nu $ is initially set to $ 1000 $, but it is gradually reduced or increased depending on the convergence of the objective function.

The algorithm starts with an initial guess, $ \bm{\tau}^{(0)} $, and the value of $ \bm{\tau} $ at the $ n$-th iteration is updated as 
\begin{equation}
	\bm{\tau}^{(n+1)} = \bm{\tau}^{(n)} + \eta^{(n)}\bm{h}^{(n)},
\end{equation}
where $ \eta^{(n)} = 1$ is used unless a negative weight is encountered. If a negative weight is encountered at the $ i $-th entry of $ \bm{\tau}^{(n+1)} $, then the update is recomputed using
\begin{equation}
	\eta^{(n)} = \frac{(\varepsilon - \bm{\tau}^{(n)}_{i})}{\bm{h}^{(n)}_{i}},
\end{equation}
where $ \varepsilon > 0$ is an arbitrary lower bound for the update of the negative weight, and is set to $ \varepsilon = 10^{-4} $ in all of our cases.

Despite being more robust than Newton's method, the LMA still suffers from bad initial guesses; especially, as the number of parameters grows and the quadrature accuracy increases, often stagnating at a suboptimal local minimum. To mitigate these issues, the LMA is coupled with a particle swarm optimization (PSO) algorithm. The PSO algorithm starts with an initialization of $ n_{c} $ particles, each with a random initial guess of $ \bm{\tau} $.  The objective function, $ \bm{g} $, given in \cref{eq:prob matrix}, is computed for each particle and the personal best, $ \bm{\tau}_{pb}  $, and global best, $ \bm{\tau}_{gb}  $, approximations are tracked throughout the iterations. The PSO step size or velocity vector, $ \bm{v} $, which is initialized as the zero vector, is computed for each particle as 
\begin{equation}
	\widetilde{\bm{v}}^{(n+1)} = b\widetilde{\bm{v}}^{(n)}+ c_1 \bm{r}_{1}\circ\left(\widetilde{\bm{\tau}}^{(n)}_{pb} - \widetilde{\bm{\tau}}^{(n)}\right)+c_2 \bm{r}_{2}\circ\left(\widetilde{\bm{\tau}}^{(n)}_{gb} - \widetilde{\bm{\tau}}^{(n)}\right),
\end{equation}
where $ b=0.6 $ is the inertial weight parameter, $ c_1 = 1.5$ is the cognitive parameter, $ c_2 = 1.5 $ is the social parameter, $ \bm{r}_1 $ and $ \bm{r}_2 $ are vectors of length equal to that of $ \widetilde{\tau} $ with uniform random entries on $ [0,1] $, and $ \circ $ denotes an elementwise multiplication. The vector, ${\bm{\tau}}$,  is updated at each iteration as
\begin{equation}
	\bm{\tau}^{(n+1)} = \bm{\tau}^{(n)} + \bm{v}^{(n+1)}.
\end{equation} 
If a negative weight is encountered for any particle, it is simply replaced by a small positive number, \eg, $ 10^{-4} $, and the update is recomputed. As the quadrature accuracy is increased, the algorithm sometimes stagnates at a local minimum and further exploration is hindered. If the same local minimum is obtained over a number of iterations, then $ \widetilde{\bm{\tau}} $ is perturbed as $ (1-\delta)\widetilde{\bm{\tau}} + \delta\bm{r}$, where $ \bm{r} $ is a vector of length equal to the length of $ \widetilde{\bm{\tau}} $ with uniform random entries on $ [0,1] $ and $ \delta>0 $ is an arbitrary small number. The choice of the value of $ \delta $ depends on the quadrature degree, as the sensitivity to perturbation increases with the quadrature degree.

The PSO and LMA are coupled in such a way that the output vector of one is used as an initial vector of the other in a loop until convergence. The PSO mitigates issues related to initial guesses and convergence to suboptimal local minima, while the LMA offers fast convergence when good initial values are provided. Despite the efficiency of the coupled algorithm, at very high quadrature degrees, the minimization sometimes stagnates before convergence is realized. In such cases, the minimization is restarted, and in some instances, the parameters associated with the interior nodes are initialized using parameters of known PI rules.

\section{Quadrature rules and SBP operators}\label{sec:quad rules}
An SBP operator on a compact reference domain, $ \hat{\Omega}$, with a piecewise smooth boundary, $ \hat{\Gamma} $, is defined as \cite{hicken2016multidimensional}:
\begin{definition} \label{def:sbp}
	$\Dxi\in \IRtwo{n_p}{n_p}$ is a degree $ p $ SBP operator in the $ i $-direction approximating the first derivative $ \pdv{\hat{\bm{x}}_i} $ on the set of nodes $ {S}=\{\bm{x}_{j}\}_{j=1}^{n_p} $ if
	\begin{enumerate}
		\item $[ \Dxi \bm{p}]_j = {\pdv{\fnc{P}}{\hat{\bm{x}}_i}}({\bm{x}_j}) $  for all $\fnc{P} \in \polyref{p} $
		\item $ \Dxi={\H}^{-1}\Qxi $, where $ {\H}$ is a symmetric positive definite matrix, and 
		\item $ \Qxi + \Qxi^{T}=  \Exi$, where  
		\[ \bm{p}^T \Exi  \bm {q} = \int_{\hat{\Gamma}} \fnc{P}\fnc{Q} \;\bm{n}_{\hat{\bm{x}}_i} \dd{\Gamma},\quad \forall \fnc{P},\fnc{Q} \in  \mathbb{P}^{r}(\hat{\Gamma}), \]
		where $ r \ge p $, $ \bm{n}_{\hat{\bm{x}}_i} $ is the $ i$-component of the outward pointing unit normal vector on $\hat{\Gamma} $, and $ \mathbb{P}^{q} $ denotes a polynomial space of degree $ q $.
	\end{enumerate}  
\end{definition} 
A diagonal-norm SBP operator has a diagonal $ \H $ matrix containing the weights of a volume quadrature rule of degree at least $ q_{v} = 2p-1 $ \cite{hicken2016multidimensional,hicken2013summation}. The existence of a sufficiently accurate positive-weight quadrature rule on $ \hat{\Omega} $ is necessary and sufficient for the existence of a degree $ p $ first-derivative diagonal-norm SBP operator \cite{hicken2016multidimensional}. The boundary operator, $ \Exi $, is also constructed using a degree $ 2p $ accurate positive-weight facet quadrature rule as \cite{fernandez2018simultaneous}
\begin{equation} \label{eq: Exi}
	\Exi = \sum_{\gamma\in\hat{\Gamma}} \Rg^T \Bg \Nxig \Rg,
\end{equation}
where $ \Nxig$ is a diagonal matrix containing the $ i $-component of the outward unit normal vector on facet $ \gamma $, $ \Rg $ is an extrapolation operator from the volume to the facet nodes, and $ \Bg $ is a diagonal matrix containing the facet quadrature weights. If the volume and facet quadrature nodes are collocated, then $ \Rg $ simply picks out function values at the facet nodes, resulting in a diagonal $ \Exi $ matrix. The collocation of the facet and volume nodes reduces the cost of element coupling via simultaneous approximation terms (SATs), especially for entropy-stable discretizations. For further discussion on construction of multidimensional SBP operators, we refer the reader to \cite{hicken2016multidimensional,fernandez2018simultaneous}. Construction of SATs for SBP discretizations of various model equations in CFD can be found in \cite{fernandez2018simultaneous,crean2018entropy,shadpey2020entropy,yan2018interior,worku2021simultaneous,worku2023entropy}.
 
\begin{figure}[!t]
	\centering
	\subfloat[][$q_v=1$, $n_p=6$]
	{\includegraphics[scale=0.27]{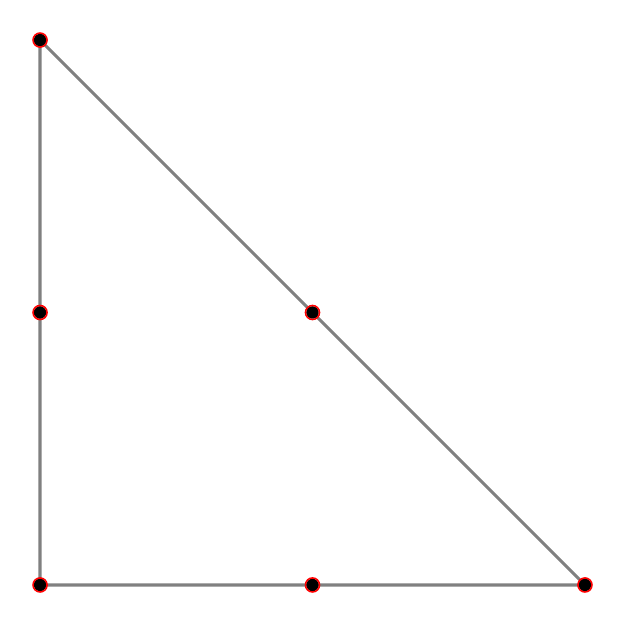} \label{fig:lgl_q1}}
	\hfill
	\subfloat[][$q_v=8$, $n_p=27$]
	{\includegraphics[scale=0.27]{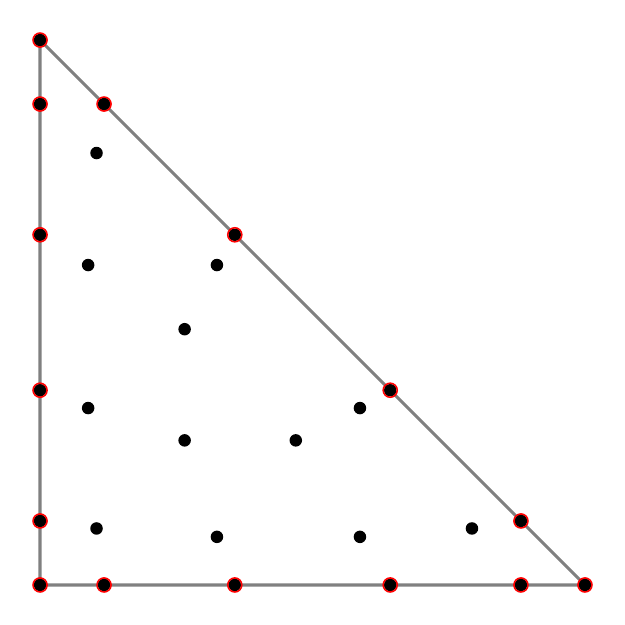}\label{fig:lgl_q8}}
	\hfill
	\subfloat[][$q_v=14$, $n_p=57$]
	{\includegraphics[scale=0.27]{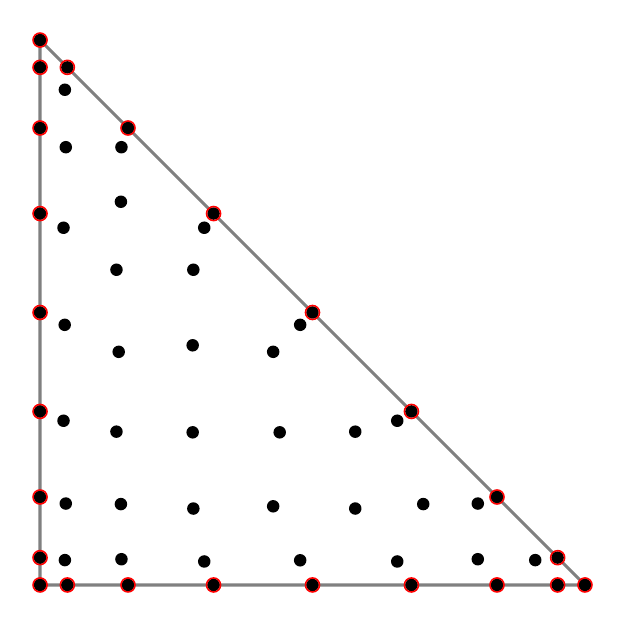}\label{fig:lgl_q14}}
	\hfill
	\subfloat[][$q_v=17$, $ n_p=78$]
	{\includegraphics[scale=0.27]{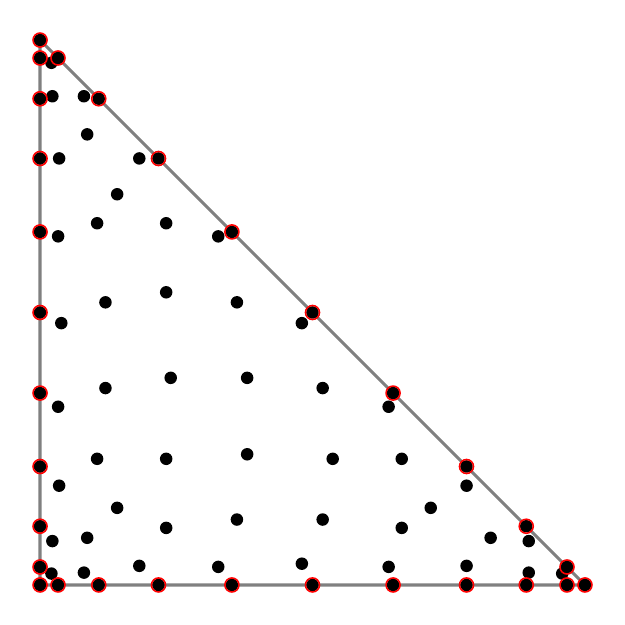} \label{fig:lgl_q17}}
	\\
	\subfloat[][$q_v=1$, $n_p=6$]
	{\includegraphics[scale=0.27]{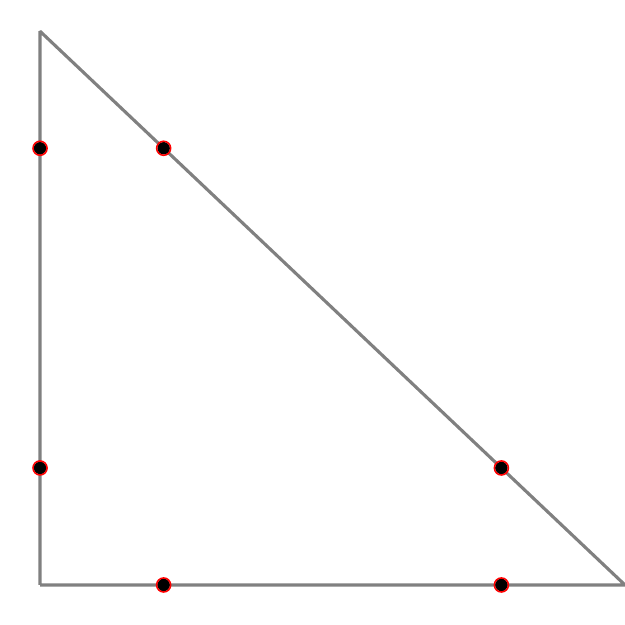}\label{fig:lg_q1}}
	\hfill
	\subfloat[][$q_v=8$, $n_p=28$]
	{\includegraphics[scale=0.27]{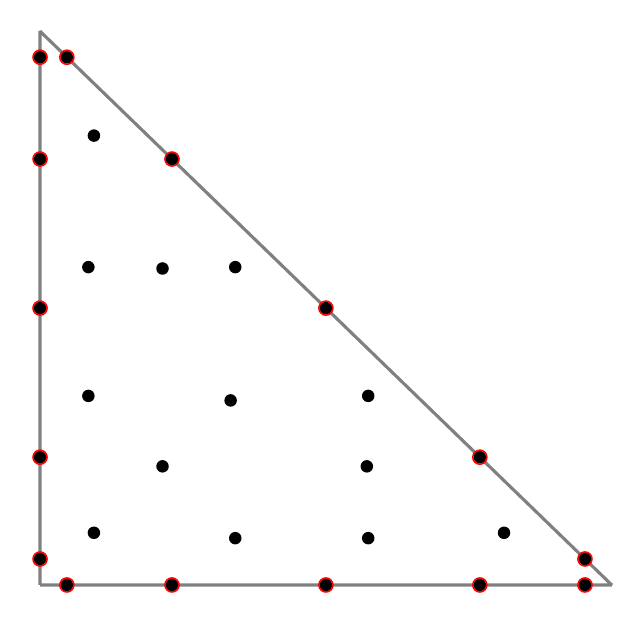}\label{fig:lg_q8}}
	\hfill
	\subfloat[][$q_v=14$, $n_p=60$]
	{\includegraphics[scale=0.27]{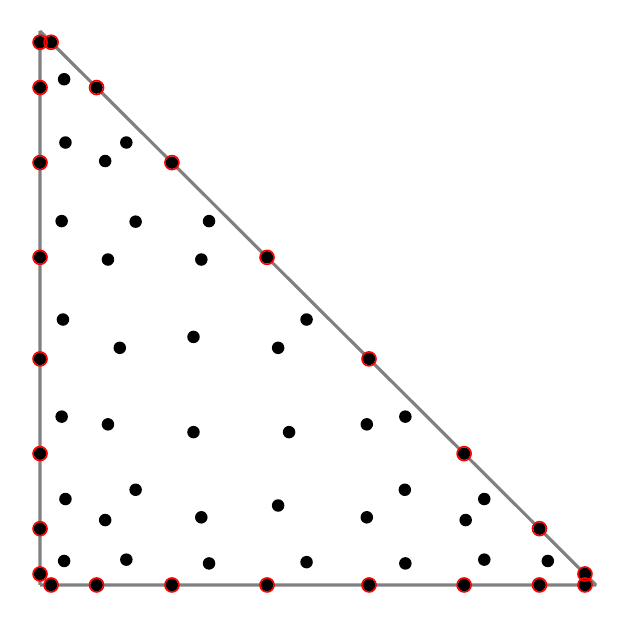} \label{fig:lg_q14}}
	\hfill
	\subfloat[][$q_v=17$, $ n_p=81$]
	{\includegraphics[scale=0.27]{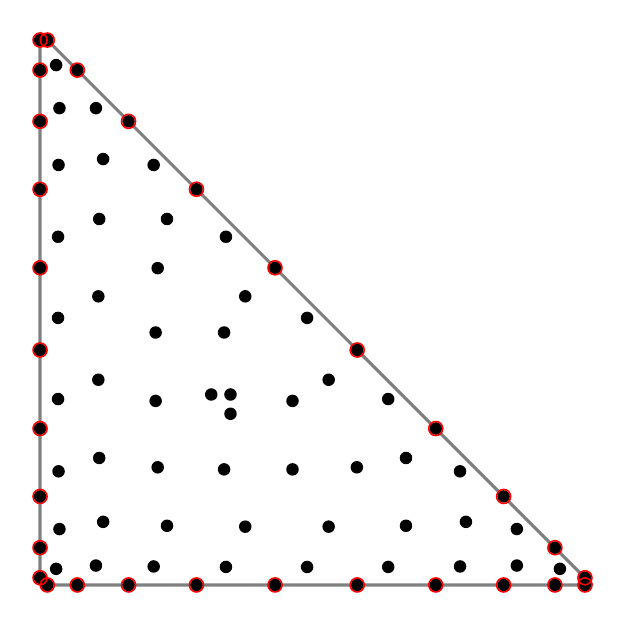} \label{fig:lg_q17}}
	\\
	\subfloat[][$q_v=1$, $n_p=6$]
	{\includegraphics[scale=0.016]{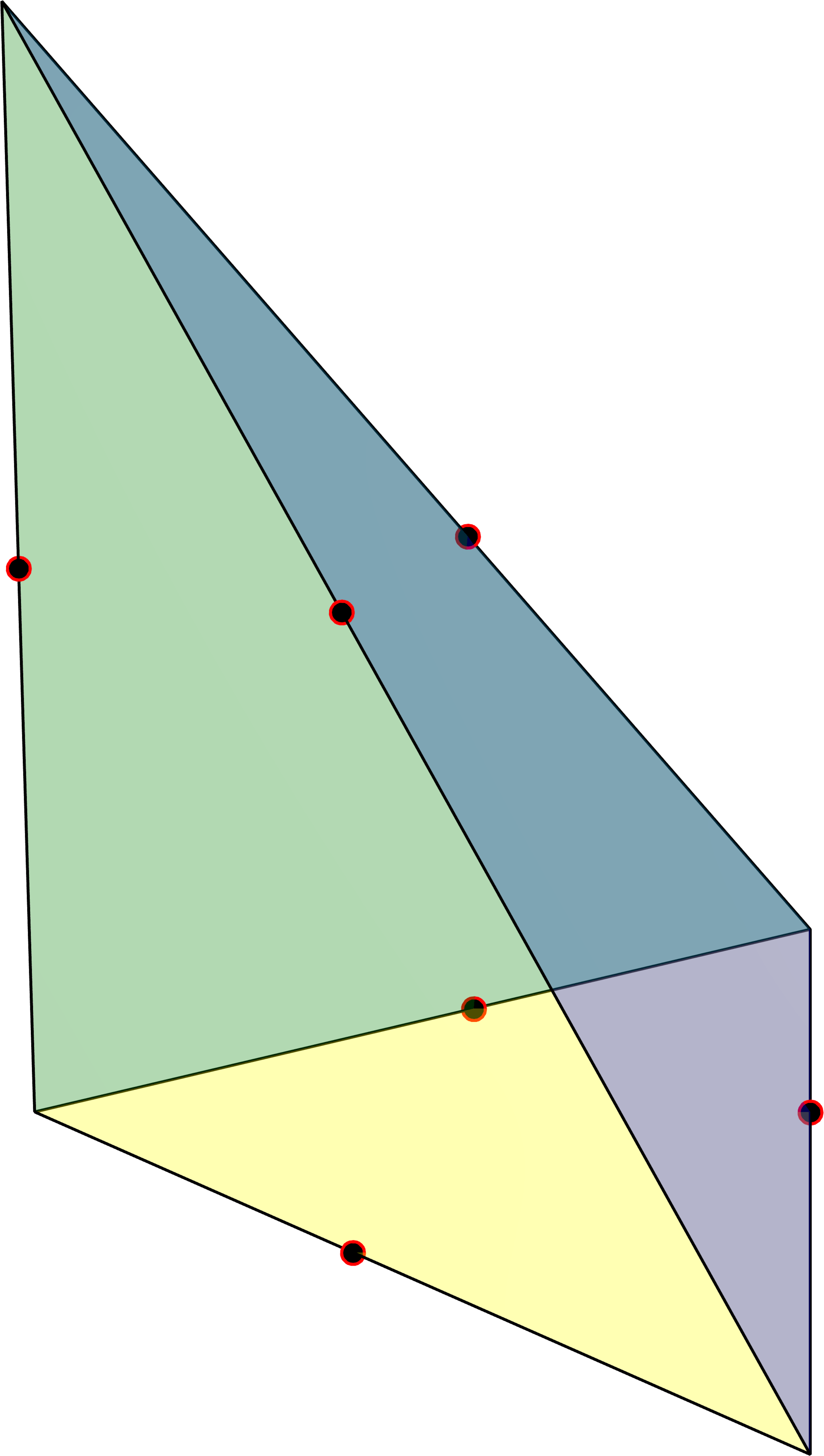}\label{fig:tet_q1}}
	\hfill
	\subfloat[][$q_v=3$, $n_p=23$]
	{\includegraphics[scale=0.016]{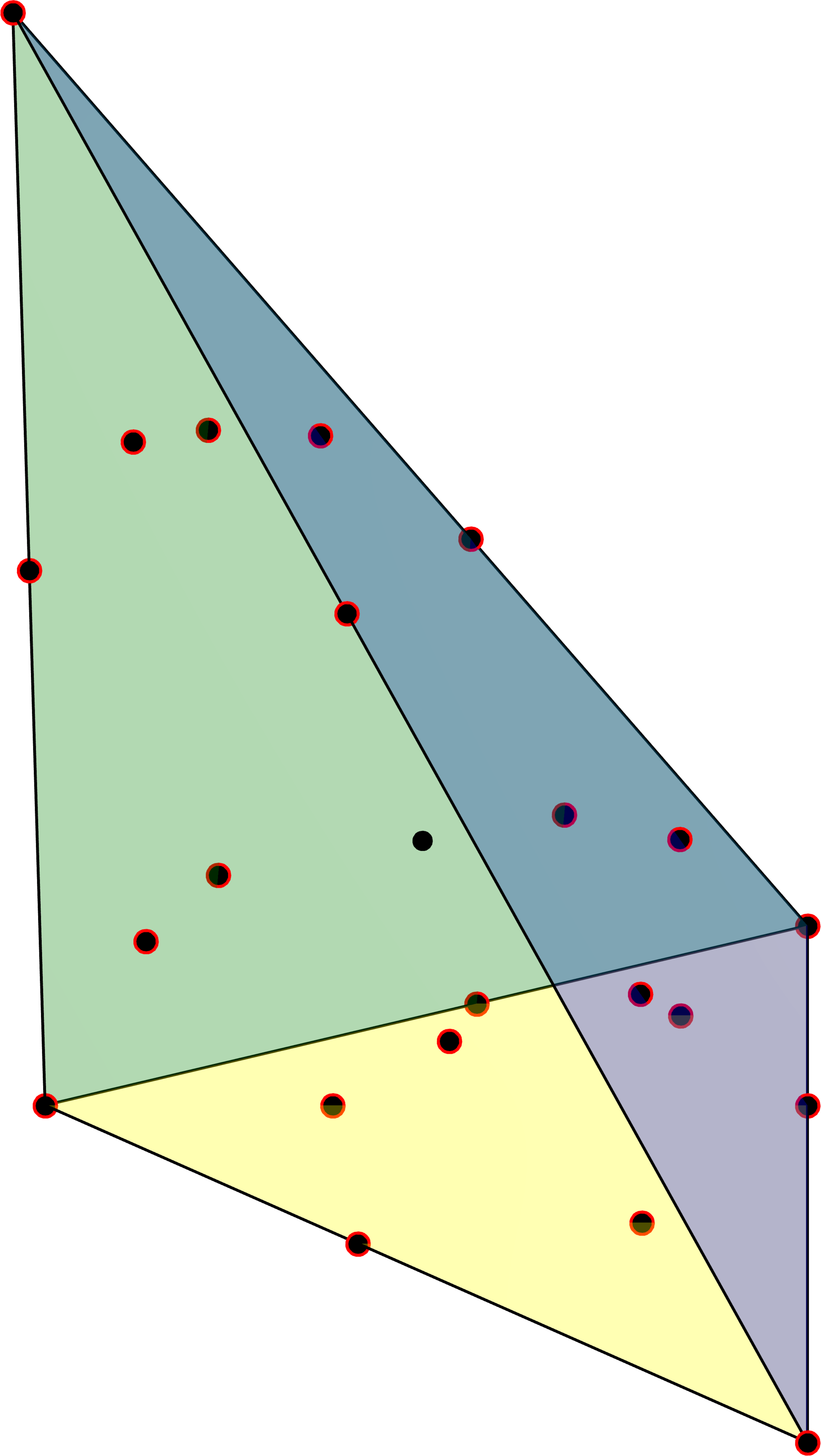} \label{fig:tet_q3}}
	\hfill
	\subfloat[][$q_v=5$, $n_p=44$]
	{\includegraphics[scale=0.016]{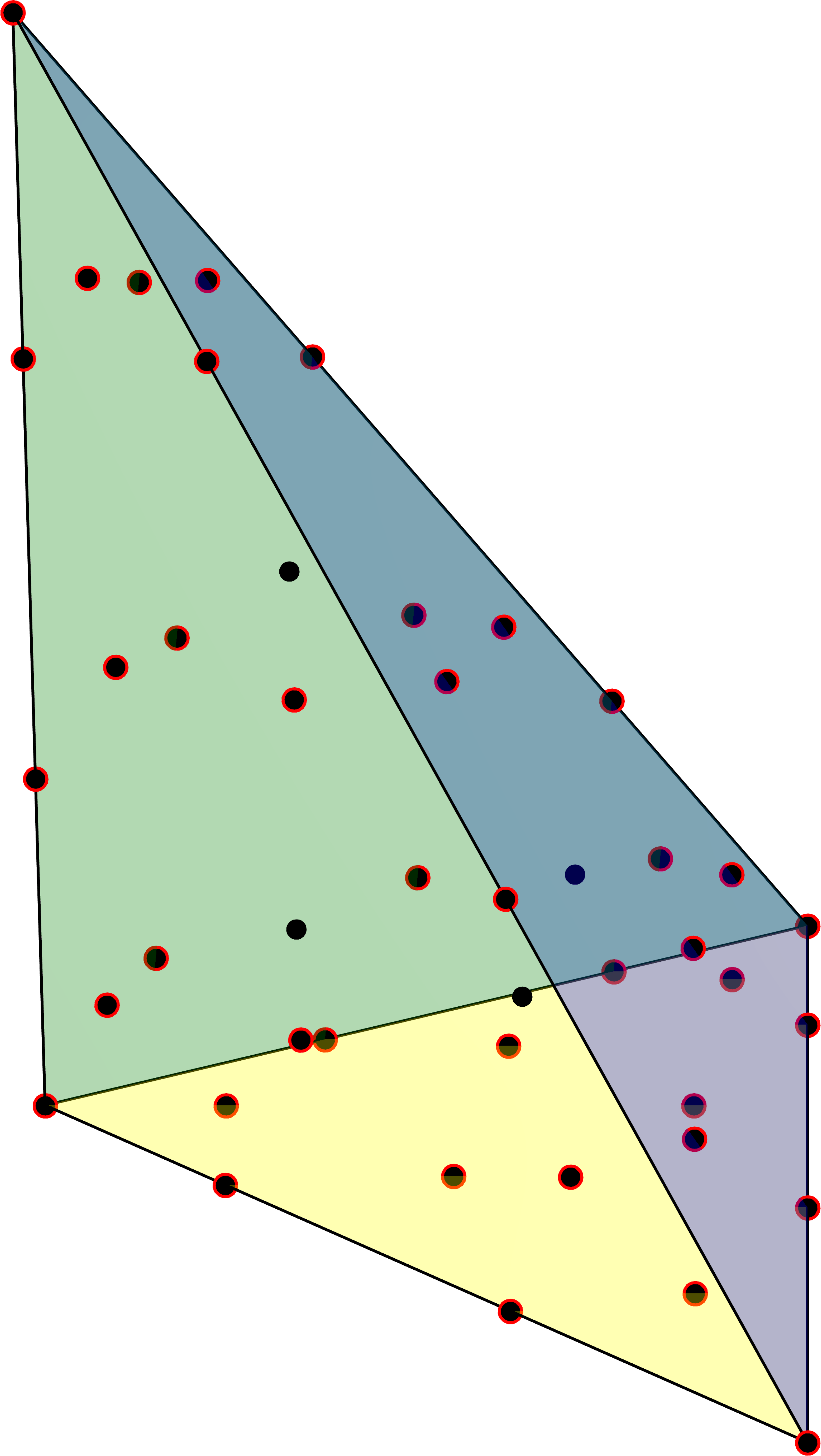}\label{fig:tet_q5}}
	\hfill
	\subfloat[][$q_v=9$, $ n_p=121$]
	{\includegraphics[scale=0.016]{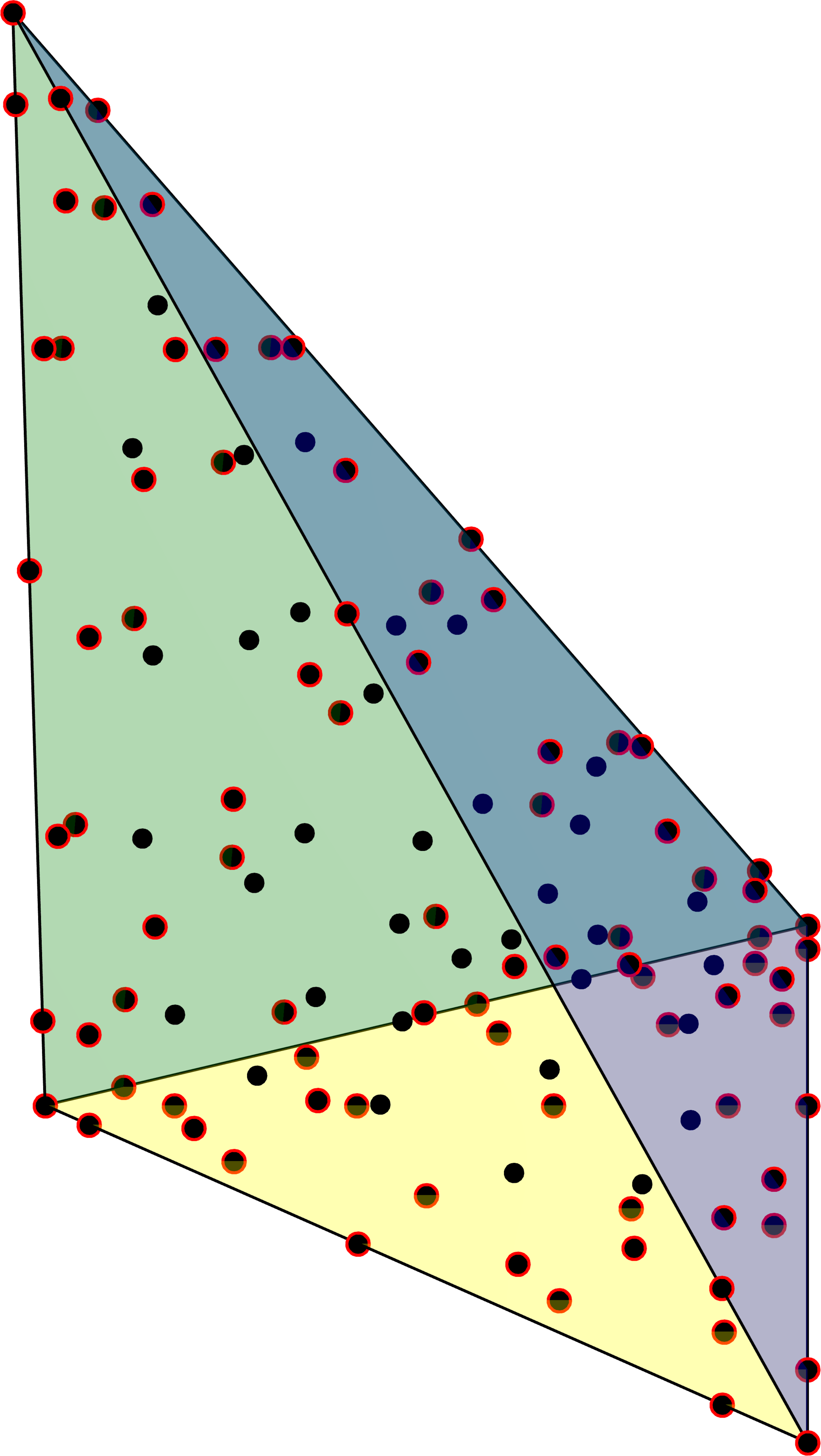} \label{fig:tet_q9}}
	\caption{\label{fig:quad_rules} Nodal locations of selected quadrature rules with LGL facet nodes (top row), LG facet nodes (middle row),  and rules on the tetrahedron (bottom row). The symbols $\sbullet$ and $\red{\bm{\circ}}$ denote the volume and facet nodes, respectively, which are collocated on the facets.}
\end{figure}

Using the methodology outlined in the previous section, we have derived quadrature rules that satisfy conditions i. -- v. stated in \cref{sec:preliminaries}, and constructed  SBP diagonal-$ \E$ operators on the reference triangle and tetrahedron. Quadrature rules with LGL and LG facet nodes of degree up to twenty on the triangle and up to ten on the tetrahedron are derived, which can be found in the supplementary data repository\footnote{\url{https://github.com/OptimalDesignLab/SummationByParts.jl/tree/master/quadrature_data}}. \cref{fig:quad_rules} illustrates the nodal configurations of some of the quadrature rules in 2D and 3D. Many of the rules are novel to the authors' knowledge, and substantial improvements, in terms of number of quadrature points, have been achieved for several of the existing rules on the tetrahedron, as illustrated in \cref{tab:number of nodes}. These improvements result in more efficient SBP diagonal-$ \E $ operators; for instance, new operators of degree $ 3$ and $ 4 $ are constructed with $ 44 $ and $ 76 $ nodes, respectively, instead of $ 69 $ and $ 99 $ nodes. The derived quadrature rules extend the available set of SBP diagonal-$ \E $ operators from degree 4 to 10 in 2D and from degree 4 to 5 in 3D. Unless specified otherwise, all the numerical studies in this work use SBP diagonal-$ \E $ operators with quadrature rules stated in the first columns of the different types of quadrature rules presented in \cref{tab:number of nodes}, which also provides the minimum nodal spacing of the quadrature rules on the reference elements. It is noted that all rules on the triangle with LGL facet nodes have larger minimum nodal spacing than those with LG facet nodes. Furthermore, the rules obtained on the tetrahedron have equal or larger minimum nodal spacing than existing rules.

\begin{table*} [!t]
	\small
	\centering
	\begin{threeparttable}
	\caption{\label{tab:number of nodes} Number of quadrature points, $ n_p $, and minimum quadrature node spacing, $ \Delta r $, on the reference triangle and tetrahedron. Quadrature rules in the cited columns denote existing results, while novel quadrature rules are underlined.}
	
	\setlength{\tabcolsep}{0.2em}
	\renewcommand*{\arraystretch}{1.2}
	\begin{tabular}{ccccccccccccccc}
		\toprule
		\multirow{2}{*}{$ q_v $} &  \multicolumn{4}{c}{Tri-LGL } &  \multicolumn{4}{c}{Tri-LG}&   \multicolumn{6}{c}{Tet} \\ 
		\cmidrule(lr){2-5} \cmidrule(lr){6-9} \cmidrule(lr){10-15}
		& $n_{p}$ & $\Delta r$& $ n_{p} $\cite{fernandez2019staggered,hicken2023summationbyparts}& $ n_{p}$\cite{wu2021high} & $n_{p}$ & $\Delta r$& $ n_{p} $\cite{hicken2023summationbyparts}& $ n_{p}$\cite{chen2017entropy} & $n_{p}$ & $\Delta r$& $ n_{p} $\cite{hicken2023summationbyparts} & $\Delta r$\cite{hicken2023summationbyparts}& $ n_{p} $\cite{marchildon2020optimization} & $\Delta r$\cite{marchildon2020optimization}\\
		\midrule
		1 & 6   &1.000	&     &  6 			& 6	   &0.423    &      & 6  		&6     &1.000    &        &    &6 & 1.000    \\
		2 & 7   &0.471	&  7 &     			& 7	    &0.423   &  7  &     		&7      &0.866    & 13  &0.333    &7  &0.866    \\
		3 & 10 &0.553	&     & 10		  &10    &0.225   &      & 10 		&    &    &       &    &23 &0.222    \\
		4 & 12 &0.325 	& 12&     		  &12	 &0.225    & 12 &      		&$ \underline{23}^{\dagger} $    &0.377    &  36 &0.183    &26& 0.222    \\
		5 & 15 &0.345	&     &15 		  &18 	 &0.083    &      & 18 		&\underline{44}    &0.237    &       &    & &    \\
		6 & 18 &0.214	& 18 &    		  &21	 &0.097   & 21 &     		&\underline{51}    &0.159    & 69   &0.106    & &    \\
		7 & 24 &0.141	&     &24 		 &22	&0.094   &      &  22 		&\underline{76}    &0.107    &       &     & &     \\
		8 & 27 &0.129	& 27&    		 &\underline{28} 	&0.094   & 28 &    			&\underline{89}    &0.094   & 99  &0.017    & &    \\
		9 & \underline{33} &0.129	&     &    		   &\underline{34}	  &0.068    &      &  			&\underline{121}    &0.077    &  &    & &    \\
		10 & \underline{36} &0.076	&    &   		   &\underline{39}	  &0.068   &       &  		   &$\underline{145}^{\ddagger}$   &0.070    &  &    & &    \\
		11 & \underline{40} &0.128	&    &    		   &\underline{42} 	  & 0.051   &       &  			 &    &    &  &    & &    \\
		12 & 48 &0.069	&48&     		&\underline{49}	   &0.051   &       &  			 &    &    &  &    & &    \\
		13 & \underline{55} &0.039	&   &     		   &\underline{54}	  &0.024   &      &  			&    &    &  &    & &    \\
		14 & \underline{57} &0.062	&    &    		   &\underline{60}	  &0.040   & 	   &  			&    &    &  &    &&     \\
		15 & \underline{69} &0.074	&    &             &\underline{69}	  &0.032    &      &  			&    &    &  &    & &    \\
		16 & \underline{72} &0.052	& 75&    		 &\underline{72} 	&0.032    & 	  &  		  &    &   &    &    & &    \\
		17 & \underline{78} &0.030	&    &    		   &\underline{81} 	   &0.026    &  	&  			&    &    &  &    & &    \\
		18 & \underline{93} &0.041	&    &    		   &\underline{93} 	  &0.026    &  		&  			&    &    &  &    &&     \\
		19 & \underline{96} &0.037	&    &    		   &\underline{96} 	  &0.010    &  		&  			&    &    &  &    & &    \\
		20 & \underline{103} &0.037&   &   		     &\underline{103} 	&0.022    &  	   &  		   &    &    &  &    & &    \\
		\bottomrule
	\end{tabular}
	\begin{tablenotes}
		\item[$ \dagger $] The $ q_v=3 $ and $ q_v=4 $ quadrature rules on the tetrahedron are identical.
		\item[$ \ddagger $]A rule with 139 nodes is provided in the supplementary data repository, but it is not considered here since it leads to a restrictively small timestep.
	\end{tablenotes}
	\end{threeparttable}
\end{table*}

\section{Numerical results}\label{sec:num results}
In this section, the SBP diagonal-$ \E $ operators constructed using the proposed quadrature rules are applied to linear and nonlinear problems. First, a mesh is generated by partitioning the spatial domain, $\Omega$, into $m$ squares or cubes in each direction and subdividing them into two triangles or six tetrahedra, respectively. The nodes on the physical elements are obtained by affinely mapping the nodes on the reference elements. The 2D and 3D meshes are refined for mesh convergence studies using $ m_{k} = 60 - 5p + (12-p)k$ and $ m_{k} = 10 + 5k $ number of edges in each direction, respectively, where $ k=\{1,2,\dots\} $ denotes the refinement level. The number of edges are chosen to ensure that errors are sufficiently larger than machine precision, enabling calculations of convergence rates for the highest-degree operators.

The standard fourth-order Runge-Kutta (RK4) scheme is applied to march the numerical solution in time. For the accuracy studies, sufficiently small timesteps are used such that the temporal errors are negligible compared to the spatial errors. As in \cite{shadpey2020entropy,worku2023entropy}, the $ L^2 $ solution error in the domain is computed by interpolating the numerical solution from the SBP nodes to a quadrature rule of degree $ 3p + 1 $, integrating the square of the solution error, summing the result over all elements, and taking the square root of the sum. 

\subsection{Linear advection problem}
We consider the linear advection equation,
\begin{equation}\label{eq:adv eq}
	\pder[\fn{U}]{t} + \sum_{i=1}^{d}\bm{c}_{i}\pder[\fn{U}]{ \hat{\bm{x}}_{i}} = 0,
\end{equation} 
on the periodic domain $ \Omega=[0,1]^{d} $. The problem is used to test the accuracy and timestep stability limits of the operators. The initial condition is obtained from the exact solution,
\begin{equation}\label{eq:exact adv sol}
	\fn{U}( \hat{\bm{x}}) =\prod_{i=1}^{d} \sin(\omega\pi (\hat{\bm{x}}_{i}-\bm{c} _{i}t)),
\end{equation}
where $ \bm{c}=[5/4,\sqrt{7}/4]^T $ in 2D or $ \bm{c}=[3/2,1/2,1/\sqrt{2}]^T $ in 3D is used in all cases. The values of $\bm{c}$ are chosen to set the wave speed magnitude at $\sqrt{d}$ but are otherwise chosen randomly. The direction of the wave propagation depends on $ \bm{c} $ and affects numerical errors and mesh convergence rates in some cases.

The advection equation, \cref{eq:adv eq}, is discretized using the diagonal-$ \E $ SBP operators and an upwind SAT, see, \eg, \cite{fernandez2018simultaneous} for the details of the SBP-SAT discretization. 
The problem is run up to $ t=1 $ with the $ \omega $ parameters in \cref{eq:exact adv sol} set to $ 8 $ and $ 2 $ for the 2D and 3D cases, respectively. The solution errors and convergence rates are tabulated in \cref{tab:conv study}, which shows convergence rates close to $ p+1 $ on the finest meshes.

\begin{table*} [!t]
	\small
	\centering
	\caption{\label{tab:conv study} Solution convergence of the linear advection problem discretized with SBP diagonal-$ \E $ operators on triangles and tetrahedra.}
	
	\setlength{\tabcolsep}{0.2em}
	\renewcommand*{\arraystretch}{1.3}
	\begin{tabular}{lllccccccccc} 
		\toprule
		\multirow{2}{*}{$ p $} &\multirow{2}{*}{$ q_v $} &  &\multicolumn{3}{c}{Tri-LGL} &  \multicolumn{3}{c}{Tri-LG}&   \multicolumn{3}{c}{Tet} \\ 
		\cmidrule(lr){4-6} \cmidrule(lr){7-9} \cmidrule(lr){10-12}
		& & &$m_1$ & $ m_2 $& $ m_3 $ & $m_1$ & $ m_2 $& $ m_3 $ & $m_1$ & $ m_2 $& $ m_3 $\\
		\midrule
		\multirow{2}{*}{1} &\multirow{2}{*}{1} & error & 3.81e-01	&  2.98e-01 &  2.37e-01  &  1.69e-01  &  1.28e-01 & 9.99e-02  &   7.74e-02 & 4.72e-02  &  3.13e-02  \\
		& & rate  & 	-- 		&  1.60 & 1.72 &	-- 		& 1.80  & 1.86 &	-- 	     & 1.72  &   1.84   \\
		\multirow{2}{*}{1} &\multirow{2}{*}{2} & error & 3.58e-02	&  2.36e-02 & 1.64e-02   &  4.15e-02  &  2.77e-02 &  1.94e-02 & 8.80e-03  & 4.42e-0  &  2.64e-03  \\
		& & rate  & 	-- 		& 2.69  & 2.72 &	-- 		& 2.62  & 2.66 &	-- 	     &   2.39 &  2.30    \\
		\multirow{2}{*}{2} &\multirow{2}{*}{3} & error & 8.60e-04	&  4.75e-04 & 2.87e-04   &  1.23e-03  & 6.76e-04  & 4.04e-04  &   &  &    \\
		& & rate  & 	-- 		& 3.85  & 3.76  &	-- 		& 3.90  & 3.85 &	     	&    	&      \\
		\multirow{2}{*}{2} &\multirow{2}{*}{4} & error & 5.27e-04	&  3.19e-04 &  2.09e-04  &  5.56e-04  & 3.37e-04  &  2.20e-04 & 9.11e-04  &  3.56e-04 &  1.75e-04  \\
		& & rate  & 	-- 		&  3.26 &  3.17&	-- 		& 3.26  &  3.17&	-- 	     &  3.26 &  3.19    \\
		\multirow{2}{*}{3} &\multirow{2}{*}{5} & error & 4.11e-05	&  2.22e-05 & 1.31e-05   &  5.24e-05  & 2.81e-05  & 1.63e-05  &  5.92e-05 & 1.94e-05  & 8.01e-06   \\
		& & rate  & 	-- 		&  3.99 & 3.98 &	-- 		& 4.05 &  4.07&	-- 	     & 3.89  &  3.96    \\
		\multirow{2}{*}{3} &\multirow{2}{*}{6} & error & 4.18e-05	& 2.26e-05  & 1.33e-05   &  4.52e-05  &  2.44e-05 & 1.43e-05  & 4.48e-05  &  1.45e-05 &  5.98e-06  \\
		& & rate  & 	-- 		& 3.99  & 3.99 &	-- 		& 4.00  & 4.00 &	-- 	     &  3.93&  3.95    \\
		\multirow{2}{*}{4} &\multirow{2}{*}{7} & error & 4.72e-06	&  2.21e-06 &  1.15e-06  & 4.02e-06  & 1.87e-06  & 9.68e-07  & 3.44e-06  & 8.43e-07  &  2.85e-07  \\
		& & rate  & 	-- 		& 4.91  & 4.93 &	-- 		& 4.94  & 4.95 &	-- 	     &  4.89 &  4.86    \\				
		\multirow{2}{*}{4} &\multirow{2}{*}{8} & error & 3.63e-06	& 1.70e-06  &  8.82e-07  & 4.35e-06   & 2.08e-06 & 1.09e-06  & 2.70e-06  &6.28e-07   &  2.08e-07  \\
		& & rate  & 	-- 		& 4.91  & 4.93 &	-- 		& 4.79 & 4.85 &	-- 	     & 5.07  & 4.95     \\ 
		\multirow{2}{*}{5} &\multirow{2}{*}{9} & error & 3.51e-07	&  1.41e-07 & 6.30e-08   & 4.31e-07   & 1.72e-07  &  7.81e-08 &  1.44e-07 &  2.64e-08 &  7.14e-09  \\
		& & rate  & 	-- 		& 5.94  &  6.01&	-- 		&  5.95 & 5.92 &	-- 	     & 5.90  & 5.86     \\
		\multirow{2}{*}{5} &\multirow{2}{*}{10}& error & 3.65e-07	&  1.49e-07 & 6.88e-08   & 3.54e-07   & 1.43e-07  & 6.51e-08  & 9.62e-08  & 1.78e-08  &  4.80e-09  \\
		& & rate  & 	-- 		&  5.82 & 5.78 &	-- 		& 5.90  & 5.89 &	-- 	     &  5.87 & 5.86     \\
		\multirow{2}{*}{6} &\multirow{2}{*}{11} & error & 4.27e-08 &  1.44e-08 & 5.52e-09   & 4.60e-08   & 1.54e-08  &  5.95e-09 &   &   &    \\
		& & rate  & 	-- 		& 7.07  & 7.16 &	-- 		& 7.08  &  7.14 &	     &   &      \\
		\multirow{2}{*}{6} &\multirow{2}{*}{12} & error & 4.68e-08 & 1.58e-08  & 6.16e-09   & 4.59e-08   & 1.57e-08  & 6.21e-09  &   &   &    \\
		& & rate  & 	-- 		&  7.05 & 7.04 &	-- 		& 6.96 & 6.94 &	 	     &   &      \\
		\multirow{2}{*}{7} &\multirow{2}{*}{13} & error & 1.26e-08 & 3.62e-09  &  1.23e-09  &  1.34e-08  & 3.83e-09  & 1.30e-09  &   &   &    \\
		& & rate  & 	-- 		& 8.08  &  8.07&	-- 		&  8.13 &  8.11&	     &   &      \\
		\multirow{2}{*}{7} &\multirow{2}{*}{14}& error & 1.25e-08	&  3.61e-09 & 1.23e-09   & 1.29e-08   & 3.77e-09  & 1.30e-09  &   &   &    \\
		& & rate  & 	-- 		&  8.06 & 8.05 &	-- 		& 8.00  & 7.99 &	     &   &      \\
		\multirow{2}{*}{8} &\multirow{2}{*}{15} & error & 4.46e-09 &  1.16e-09 & 3.57e-10   & 4.55e-09   &  1.20e-09 & 3.74e-10  &   &   &    \\
		& & rate  & 	-- 		& 8.72  & 8.85 &	-- 		&  8.65 &  8.72&	     &   &      \\
		\multirow{2}{*}{8} &\multirow{2}{*}{16} & error & 4.93e-09 & 1.30e-09  & 4.04e-10   &  4.80e-09  & 1.24e-09  & 3.80e-10  &   &   &    \\
		& & rate  & 	-- 		& 8.65  & 8.75 &	-- 		&  8.77 & 8.87 &	     &   &      \\
		\multirow{2}{*}{9} &\multirow{2}{*}{17} & error & 5.48e-09 & 1.36e-09  &  3.63e-10  &  4.91e-09  &  1.17e-09 & 3.22e-10  &   &   &    \\
		& & rate  & 	-- 		& 9.06  &  9.87 &	-- 		& 9.28  &  9.69& 	     &   &      \\
		\multirow{2}{*}{9} &\multirow{2}{*}{18} & error & 4.78e-09 & 1.14e-09  &  3.07e-10  &  5.02e-09  &  1.21e-09 & 3.32e-10  &   &   &    \\
		& & rate  & 	-- 		&  9.32 & 9.79 &	-- 		&  9.22 & 9.70 & 	     &   &      \\
		\multirow{2}{*}{10} &\multirow{2}{*}{19} & error & 	2.59e-08 &  5.03e-09 & 1.19e-09   &  2.73e-08  &  5.29e-09 & 1.25e-09  &   &   &    \\
		& & rate  & 	-- 		&  10.63 & 10.77 &	-- 		&  10.64 & 10.82 & 	     &   &      \\
		\multirow{2}{*}{10} &\multirow{2}{*}{20}& error & 2.63e-08 & 5.17e-09  & 1.24e-09   &  2.62e-08  & 5.10e-09  & 1.21e-09  &   &   &    \\
		& & rate  & 	-- 		& 10.55  & 10.71 &	-- 		&  10.62 & 10.76 & 	     &   &      \\
		\bottomrule
	\end{tabular}
\end{table*}

In addition to the accuracy of the operators, we are also interested in their maximum timestep limits for explicit time marching schemes. A large timestep limit is desired for stability-bounded problems, where the maximum stable timestep can be applied without compromising accuracy. The maximum timestep is computed for each operator using golden section optimization. For this study, the triangular and tetrahedral meshes are obtained by subdividing quadrilateral and hexahedral meshes with four elements in each direction. The discretization is considered to be stable if the change in energy is less than or equal to zero after five periods. The change in energy at a given timestep is computed as,
\begin{equation}
	\Delta E = \sum_{\Omega_{k}\in \fn{T}_h} \left(\bm{u}^T_k\H_k\bm{u}_k - \bm{u}_{0,k}^T\H_{k}\bm{u}_{0,k}\right), 
\end{equation}
where $ \bm{u}_{k} $, $ \bm{u}_{0,k}$, and $ \H_k $ are the solution vector at the specified timestep, the initial solution vector, and the norm matrix on element $ \Omega_{k} $, respectively, and $ \fn{T}_h $ is a set containing all physical elements. 

\cref{tab:timestep adv} presents the maximum timestep values for each SBP-SAT discretization. On the triangle, we have not made improvements on existing quadrature rules except in the case of the degree $ 16 $ rule with the LGL facet nodes; hence, comparisons of the maximum timesteps with previously existing operators are not presented. On the tetrahedron, the new quadrature rule for the degree 2 SBP diagonal-$\E$ operator yields a smaller stable timestep than the existing rules, while the degree 3 and 4 operators with $ q_v=2p $ quadrature rules lead to slightly lower but comparable stable timesteps relative to the existing rules. The degree 3 and 4 diagonal-$ \E $ operators constructed with $ q_v=2p-1 $ quadrature rules have about 1.72 and 1.28 times larger maximum timesteps, respectively, than the existing degree 3 and 4 operators. This, combined with their lower node count, leads to substantial efficiency improvements for stability-bounded problems.

\begin{table*} [!t]
	\small
	\centering
	\caption{\label{tab:timestep adv} Maximum timestep on the reference triangle and tetrahedron.}
	\begin{threeparttable}
	\setlength{\tabcolsep}{0.75em}
	\renewcommand*{\arraystretch}{1.2}
	\begin{tabular}{cccccc}
		\toprule
		\multirow{1}{*}{$ q_v $} &  \multicolumn{1}{c}{Tri-LGL } &  \multicolumn{1}{c}{Tri-LG}&   \multicolumn{1}{c}{Tet}  &   \multicolumn{1}{c}{Tet\cite{hicken2023summationbyparts}} &  \multicolumn{1}{c}{Tet\cite{marchildon2020optimization}}\\ 
		
		\midrule
		1    & 0.0258   & 0.0758   &0.0641   &     			 &0.0641    \\
		2   & 0.0446   & 0.0394   & 0.0388  & 0.0345 & 0.0388  \\
		3   & 0.0323   & 0.0365   &    			 &     			  & 0.0164   \\
		4   & 0.0235   & 0.0216   &0.0065  &  0.0208  & 0.0204   \\
		5   & 0.0258   & 0.0203   & 0.0184  &     			 &    \\
		6   &0.0146    & 0.0131   & 0.0097  & 0.0101    &    \\
		7   &  0.0145   &  0.0129  &0.0074   &     			&    \\
		8   & 0.0106   & 0.0092   & 0.0049  & 0.0058 &    \\
		9    & 0.0070    &  0.0089  & 0.0075  &     	    &    \\
		10   &0.0066    &  0.0012   & 0.0059  &     		&    \\
		11    & 0.0008   & 0.0066   &	&   &	 \\
		12   &0.0059    & 0.0015    &	&   &	 \\
		13   &0.0016    &0.0040    &	&   &	 \\
		14   &0.0048    & 0.0040   &	&   &	\\
		15   & 0.0025   & 0.0035   &	&   &	\\
		16   & 0.0038   & 0.0020   &	&   &	 \\
		17    & 0.0022   & 0.0027   &	&   &	 \\
		18   &0.0033    & 0.0024   &	&   &	\\
		19   &0.0029    & 0.0023   &	&   &	\\
		20   &0.0025    & 0.0019   &	&   &	 \\
		\bottomrule
	\end{tabular}
\end{threeparttable}
\end{table*}

\subsection{Isentropic vortex problem}
The isentropic vortex problem, governed by the Euler equations, is another common test case used to study the accuracy of high-order methods. We consider the 3D case on the periodic domain $\Omega = \left[-10,10\right]^{3}$ with the initial conditions \cite{williams2013nodal}
\begin{equation}
	\begin{aligned}
		\rho & =\left(1-\frac{2}{25}\left(\gamma-1\right)\exp\left(1-\left(x_2-t\right)^{2}-x_1^{2}\right)\right)^{\frac{1}{\gamma-1}}, \\
		e & =\frac{\rho^{\gamma}}{\gamma(\gamma-1)}+\frac{\rho}{2}\left(u^{2}+v^{2}+w^{2}\right), \\
		u & =-\frac{2}{5}\left(x_2-t\right)\exp\left(\frac{1}{2}\left[1-\left(x_2-t\right)^{2}-x_1^{2}\right]\right), \\ 
		v & =1+\frac{2}{5}x_1\exp\left(\frac{1}{2}\left[1-\left(x_2-t\right)^{2}-x_1^{2}\right]\right), \\
		w & =0,
	\end{aligned}
\end{equation}
where $ \rho $ is the density, $ e $ is the total energy per unit volume, $ u $, $ v $, and $ w $ are the velocities in the $x_1 $, $x_2 $, and $ x_3 $ directions, respectively, and $\gamma=7/5$ is the ratio of specific heats. 

We use the Hadamard-form entropy-stable discretization \cite{crean2018entropy} on tetrahedral elements with the Ismail-Roe two-point fluxes \cite{ismail2009affordable}. Furthermore, the matrix-type interface dissipation operator of \cite{ismail2009affordable} is applied. The problem is run until $ t=1 $, and a mesh convergence study is conducted. The $ L^2 $ solution errors and their rates of convergence are shown in \cref{tab:conv study vortex}. Convergence rates greater than $ p+0.5 $ are attained on the finest meshes for all operators constructed using the new quadrature rules.

\begin{table*} [!t]
	\small
	\centering
	\caption{\label{tab:conv study vortex} Solution convergence of the 3D isentropic vortex problem discretized with SBP diagonal-$ \E $ operators on tetrahedra.}
	
	\setlength{\tabcolsep}{0.4em}
	\renewcommand*{\arraystretch}{1.2}
	\begin{tabular}{lllccccc}
		\toprule
		\multirow{1}{*}{$ p $} &\multirow{1}{*}{$ q_v $} &  & $m_1$ &$m_2$ & $ m_3 $& $ m_4 $ & $m_5$ \\
		\midrule
		\multirow{2}{*}{1} &\multirow{2}{*}{1} & error &  1.53e+00 &9.93e-01	& 7.03e-01	&5.27e-01	& 4.14e-01\\
		& & rate  & 	-- 	&1.51	&1.55&1.57	&1.58  \\
		\multirow{2}{*}{1} &\multirow{2}{*}{2} & error & 1.35e+00 &7.98e-01	& 5.65e-01 &4.12e-01 & 3.16e-01		\\
		& & rate  & 	-- 	& 1.82	&1.55	&1.73 &1.72	 \\
		\multirow{2}{*}{2} &\multirow{2}{*}{4} & error & 3.73e-01 &1.93e-01	& 1.06e-01	&6.68e-02	& 4.52e-02	\\
		& & rate  & 	-- 	&2.30	&2.68	&2.53	&2.54	 \\
		\multirow{2}{*}{3} &\multirow{2}{*}{5} & error & 9.12e-02 &3.59e-02	& 1.68e-02	&8.87e-03	& 5.03e-03	\\
		& & rate  & 	-- 	& 3.24	&3.40 &3.51	&3.68 \\
		\multirow{2}{*}{3} &\multirow{2}{*}{6} & error & 9.27e-02 &3.68e-02& 1.72e-02&8.89e-03	& 5.02e-03	\\
		& & rate  & 	-- 	&3.21	&3.40	&3.63	&3.71	\\
		\multirow{2}{*}{4} &\multirow{2}{*}{7} & error & 2.72e-02 &7.87e-03	& 2.70e-03	&1.22e-03	& 5.80e-04	\\
		& & rate  & 	-- 	&4.31&4.80 &4.36	&4.81\\					
		\multirow{2}{*}{4} &\multirow{2}{*}{8} & error & 2.82e-02 &8.25e-03	& 2.87e-03	&1.26e-03	& 6.02e-04\\
		& & rate  & 	-- 	&4.27&4.72	&4.53	&4.78\\ 
		\multirow{2}{*}{5} &\multirow{2}{*}{9} & error &6.03e-03  &1.35e-03	& 4.27e-04	&1.51e-04	& 6.30e-05\\
		& & rate  & 	-- 	&5.20	&5.17	&5.71	&5.65\\
		\multirow{2}{*}{5} &\multirow{2}{*}{10}& error &6.47e-03  &1.47e-03	& 4.57e-04& 1.60e-4& 6.63e-05\\
		& & rate  & 	-- 	&5.15&5.25	&5.77	&5.69\\
		\bottomrule
	\end{tabular}
\end{table*}

\section{Conclusions}\label{sec:conclusions}
Several novel quadrature rules that are applicable for construction of diagonal-norm diagonal-$ \E $ SBP operators on triangles and tetrahedra are derived. The quadrature rules are obtained by coupling the LMA and PSO methods to solve the nonlinear systems of equations arising from the quadrature accuracy conditions. The LMA provides fast convergence when a good initial condition is provided, while the PSO enables efficient exploration of the design space while also mitigating stagnation issues at suboptimal local minima. Quadrature rules of degrees one through twenty on triangles with both the LGL and LG type facet nodes are presented. For tetrahedra, quadrature rules of degree one through ten are reported, which, in most cases, have substantially fewer nodes than previously known rules for SBP diagonal-$ \E $ operators. The newly derived quadrature rules lead to SBP diagonal-$ \E $ operators with enhanced efficiency relative to many of the existing operators of the same degree. They also extend the available set of SBP diagonal-$ \E $ operators from degree 4 to 10 in 2D and from degree 4 to 5 in 3D.

The diagonal-norm diagonal-$ \E $ multidimensional SBP operators are applied to solve the linear advection and isentropic vortex problems on periodic domains. Mesh refinement studies for the problems show that convergence rates on the finest meshes are greater than $ p+0.5 $ in all cases.

We have found that the quadrature points with the LGL facet nodes on triangles have larger minimum nodal spacing than those with the LG facet nodes. For tetrahedra, the rules constructed in this work provide equal or larger minimum nodal spacing than those found in the literature. We have also investigated the maximum timestep values for stability, and found that most of the new rules offer comparable or larger stable timesteps relative to previously reported rules.

\section*{Declarations}
\subsection*{Conflicts of Interest}
The authors have no conflicts of interest to declare.
\subsection*{Data Availability}
All the quadrature rules reported in this work can be found in the supplementary data repository at \url{https://github.com/OptimalDesignLab/SummationByParts.jl/tree/master/quadrature_data}. The software implementation used to obtain the quadrature rules is publicly accessible at \url{https://github.com/OptimalDesignLab/SummationByParts.jl}.

\subsection*{Acknowledgments}
The authors would like to thank Professor Masayuki Yano and his Aerospace Computational Engineering Lab at the University of Toronto for the use of their software, the Automated PDE Solver (APS).  The first author would also like to thank Andr{\'e} Marchildon for the discussions on placement of quadrature points on tetrahedral elements. Computation were performed on the Niagara supercomputer at the SciNet HPC Consortium \cite{ponce2019deploying}. SciNet is funded by: the Canada Foundation for Innovation; the Government of Ontario; Ontario Research Fund - Research Excellence; and the University of Toronto. 


\bibliographystyle{spmpsci}      
\bibliography{./references}

\begin{thebibliography}{10}
\providecommand{\url}[1]{{#1}}
\providecommand{\urlprefix}{URL }
\expandafter\ifx\csname urlstyle\endcsname\relax
  \providecommand{\doi}[1]{DOI~\discretionary{}{}{}#1}\else
  \providecommand{\doi}{DOI~\discretionary{}{}{}\begingroup
  \urlstyle{rm}\Url}\fi

\bibitem{chan2018discretely}
Chan, J.: On discretely entropy conservative and entropy stable discontinuous
  {Galerkin} methods.
\newblock Journal of Computational Physics \textbf{362}, 346--374 (2018)

\bibitem{chen2017entropy}
Chen, T., Shu, C.W.: Entropy stable high order discontinuous {G}alerkin methods
  with suitable quadrature rules for hyperbolic conservation laws.
\newblock Journal of Computational Physics \textbf{345}, 427--461 (2017)

\bibitem{cohen2001higher}
Cohen, G., Joly, P., Roberts, J.E., Tordjman, N.: Higher order triangular
  finite elements with mass lumping for the wave equation.
\newblock SIAM Journal on Numerical Analysis \textbf{38}(6), 2047--2078 (2001)

\bibitem{crean2018entropy}
Crean, J., Hicken, J.E., Del Rey~Fern{\'a}ndez, D.C., Zingg, D.W., Carpenter,
  M.H.: Entropy-stable summation-by-parts discretization of the {E}uler
  equations on general curved elements.
\newblock Journal of Computational Physics \textbf{356}, 410--438 (2018)

\bibitem{fernandez2019staggered}
Del Rey~Fern{\'a}ndez, D.C., Crean, J., Carpenter, M.H., Hicken, J.E.:
  Staggered-grid entropy-stable multidimensional summation-by-parts
  discretizations on curvilinear coordinates.
\newblock Journal of Computational Physics \textbf{392}, 161--186 (2019)

\bibitem{fernandez2018simultaneous}
Del Rey~Fern{\'a}ndez, D.C., Hicken, J.E., Zingg, D.W.: Simultaneous
  approximation terms for multi-dimensional summation-by-parts operators.
\newblock Journal of Scientific Computing \textbf{75}(1), 83--110 (2018)

\bibitem{dubiner1991spectral}
Dubiner, M.: Spectral methods on triangles and other domains.
\newblock Journal of Scientific Computing \textbf{6}, 345--390 (1991)

\bibitem{felippa2004compendium}
Felippa, C.A.: A compendium of {FEM} integration formulas for symbolic work.
\newblock Engineering Computations \textbf{21}(8), 867--890 (2004)

\bibitem{fisher2013high}
Fisher, T.C., Carpenter, M.H.: High-order entropy stable finite difference
  schemes for nonlinear conservation laws: Finite domains.
\newblock Journal of Computational Physics \textbf{252}, 518--557 (2013)

\bibitem{gassner2018br1}
Gassner, G.J., Winters, A.R., Hindenlang, F.J., Kopriva, D.A.: The {BR1} scheme
  is stable for the compressible {Navier--Stokes} equations.
\newblock Journal of Scientific Computing \textbf{77}(1), 154--200 (2018)

\bibitem{giraldo2006diagonal}
Giraldo, F.X., Taylor, M.A.: A diagonal-mass-matrix triangular-spectral-element
  method based on cubature points.
\newblock Journal of Engineering Mathematics \textbf{56}, 307--322 (2006)

\bibitem{hicken2023summationbyparts}
Hicken, J.E.: {SummationByParts.jl}. {V}ersion 0.1.2 (2015).
\newblock
  \urlprefix\url{https://github.com/OptimalDesignLab/SummationByParts.jl}

\bibitem{hicken2016multidimensional}
Hicken, J.E., Del Rey~Fern{\'a}ndez, D.C., Zingg, D.W.: Multidimensional
  summation-by-parts operators: General theory and application to simplex
  elements.
\newblock SIAM Journal on Scientific Computing \textbf{38}(4), A1935--A1958
  (2016)

\bibitem{hicken2013summation}
Hicken, J.E., Zingg, D.W.: Summation-by-parts operators and high-order
  quadrature.
\newblock Journal of Computational and Applied Mathematics \textbf{237}(1),
  111--125 (2013)

\bibitem{ismail2009affordable}
Ismail, F., Roe, P.L.: Affordable, entropy-consistent {Euler} flux functions
  {II}: {Entropy} production at shocks.
\newblock Journal of Computational Physics \textbf{228}(15), 5410--5436 (2009)

\bibitem{kennedy1995particle}
Kennedy, J., Eberhart, R.: Particle swarm optimization.
\newblock In: Proceedings of ICNN'95 - International Conference on Neural
  Networks, vol.~4, pp. 1942--1948. IEEE (1995)

\bibitem{koornwinder1975two}
Koornwinder, T.: Two-variable analogues of the classical orthogonal
  polynomials.
\newblock In: R.A. Askey (ed.) Theory and Application of Special Functions, pp.
  435--495. Academic Press (1975)

\bibitem{levenberg1944method}
Levenberg, K.: A method for the solution of certain non-linear problems in
  least squares.
\newblock Quarterly of Applied Mathematics \textbf{2}(2), 164--168 (1944)

\bibitem{liu2017higher}
Liu, Y., Teng, J., Xu, T., Badal, J.: Higher-order triangular spectral element
  method with optimized cubature points for seismic wavefield modeling.
\newblock Journal of Computational Physics \textbf{336}, 458--480 (2017)

\bibitem{marchildon2020optimization}
Marchildon, A.L., Zingg, D.W.: Optimization of multidimensional diagonal-norm
  summation-by-parts operators on simplices.
\newblock Journal of Computational Physics \textbf{411}, 109380 (2020)

\bibitem{marquardt1963algorithm}
Marquardt, D.W.: An algorithm for least-squares estimation of nonlinear
  parameters.
\newblock Journal of the Society for Industrial and Applied Mathematics
  \textbf{11}(2), 431--441 (1963)

\bibitem{montoya2023efficient}
Montoya, T., Zingg, D.W.: Efficient tensor-product spectral-element operators
  with the summation-by-parts property on curved triangles and tetrahedra.
\newblock arXiv preprint arXiv:2306.05975  (2023)

\bibitem{mulder2001higher}
Mulder, W.A.: Higher-order mass-lumped finite elements for the wave equation.
\newblock Journal of Computational Acoustics \textbf{9}(02), 671--680 (2001)

\bibitem{parsani2015entropy}
Parsani, M., Carpenter, M.H., Nielsen, E.J.: {Entropy} stable discontinuous
  interfaces coupling for the three-dimensional compressible {Navier-Stokes}
  equations.
\newblock J. Comput. Physics \textbf{290}, 132--138 (2015)

\bibitem{ponce2019deploying}
Ponce, M., van Zon, R., Northrup, S., Gruner, D., Chen, J., Ertinaz, F.,
  Fedoseev, A., Groer, L., Mao, F., Mundim, B.C., et~al.: Deploying a top-100
  supercomputer for large parallel workloads: The {Niagara} supercomputer.
\newblock In: Proceedings of the Practice and Experience in Advanced Research
  Computing on Rise of the Machines (learning), pp. 1--8 (2019)

\bibitem{proriol1957family}
Proriol, J.: Sur une famille de polynomes {\`a} deux variables orthogonaux dans
  un triangle.
\newblock Comptes Rendus Hebdomadaires des S{\'e}ances de l'Acad{\'e}mie des
  Sciences \textbf{245}(26), 2459--2461 (1957)

\bibitem{ranocha2018comparison}
Ranocha, H.: Comparison of some entropy conservative numerical fluxes for the
  {Euler} equations.
\newblock Journal of Scientific Computing \textbf{76}(1), 216--242 (2018)

\bibitem{shadpey2020entropy}
Shadpey, S., Zingg, D.W.: Entropy-stable multidimensional summation-by-parts
  discretizations on hp-adaptive curvilinear grids for hyperbolic conservation
  laws.
\newblock Journal of Scientific Computing \textbf{82}(3), 70 (2020)

\bibitem{taylor2000algorithm}
Taylor, M.A., Wingate, B.A., Vincent, R.E.: An algorithm for computing {Fekete}
  points in the triangle.
\newblock SIAM Journal on Numerical Analysis \textbf{38}(5), 1707--1720 (2000)

\bibitem{williams2013nodal}
Williams, D.M., Jameson, A.: Nodal points and the nonlinear stability of
  high-order methods for unsteady flow problems on tetrahedral meshes.
\newblock In: 21st AIAA Computational Fluid Dynamics Conference, AIAA 2013-2830
  (2013)

\bibitem{wingate2008performance}
Wingate, B.A., Taylor, M.A.: Performance of numerically computed quadrature
  points.
\newblock Applied Numerical Mathematics \textbf{58}(7), 1030--1041 (2008)

\bibitem{witherden2015identification}
Witherden, F., Vincent, P.: On the identification of symmetric quadrature rules
  for finite element methods.
\newblock Computers \& Mathematics with Applications \textbf{69}(10),
  1232--1241 (2015)

\bibitem{witherden2014analysis}
Witherden, F.D., Vincent, P.E.: An analysis of solution point coordinates for
  flux reconstruction schemes on triangular elements.
\newblock Journal of Scientific Computing \textbf{61}, 398--423 (2014)

\bibitem{worku2021simultaneous}
Worku, Z.A., Zingg, D.W.: Simultaneous approximation terms and functional
  accuracy for diffusion problems discretized with multidimensional
  summation-by-parts operators.
\newblock Journal of Computational Physics \textbf{445}, 110634 (2021)

\bibitem{worku2023entropy}
Worku, Z.A., Zingg, D.W.: Entropy-split multidimensional summation-by-parts
  discretization of the {Euler and Navier-Stokes} equations.
\newblock arXiv preprint arXiv:2305.07181  (2023)

\bibitem{wu2021high}
Wu, X., Kubatko, E.J., Chan, J.: High-order entropy stable discontinuous
  {Galerkin} methods for the shallow water equations: curved triangular meshes
  and {GPU} acceleration.
\newblock Computers \& Mathematics with Applications \textbf{82}, 179--199
  (2021)

\bibitem{yan2018interior}
Yan, J., Crean, J., Hicken, J.E.: Interior penalties for summation-by-parts
  discretizations of linear second-order differential equations.
\newblock Journal of Scientific Computing \textbf{75}(3), 1385--1414 (2018)

\bibitem{zhang2009set}
Zhang, L., Cui, T., Liu, H.: A set of symmetric quadrature rules on triangles
  and tetrahedra.
\newblock Journal of Computational Mathematics pp. 89--96 (2009)

\end{thebibliography}
\addcontentsline{toc}{section}{\refname}

\end{document}